\documentclass[11pt]{amsart}

\usepackage{amsmath}
\usepackage{amssymb}
\usepackage[latin1]{inputenc} 
\usepackage{graphicx}
\usepackage{float}
\usepackage{gastex}
\usepackage[all]{xy}

\newcommand{\n}{{\noindent}}
\renewcommand{\b}{{\bigskip}}
\newcommand{\q}{{\quad}}
\renewcommand{\l}{\left}
\renewcommand{\r}{\right}
\newcommand{\R}{\Rightarrow}
\newcommand{\LR}{\Leftrightarrow}
\newcommand{\tgauche}[1]{{\vphantom{#1}}^{\mathit t}{#1}}

\newtheorem{thm}{Theorem}[]

\newtheorem{lem}{Lemma}[]
\newtheorem{rem}{Remark}[]
\newtheorem{prop}{Proposition}[]
\newtheorem{cor}{Corollary}[]
\theoremstyle{definition}

\title{Codes and noncommutative stochastic matrices}

\author[S.~Lavallée]{Sylvain Lavallée}
\address{Sylvain Lavallée:
Département de mathématiques, Université du Québec à Montréal}
\email{lavallee.sylvain.2@courrier.uqam.ca}

\author[C.~Reutenauer]{Christophe Reutenauer}
\address{Christophe Reutenauer:
Département de mathématiques, Université du Québec à Montréal}
\email{reutenauer.christophe@uqam.ca}

\author[V.~Retakh]{Vladimir Retakh}
\address{Vladimir Retakh: Department of Mathematics, Rutgers University, Piscataway, NJ 08854, USA}
\email{vretakh@math.rutgers.edu}

\author[D.~Perrin]{Dominique Perrin}
\address{Dominique Perrin: Institut Gaspard-Monge, Université de Marne-la-Vallée, France}
\email{perrin@univ-mlv.fr}

\begin{document}

\begin{abstract}
Given a matrix over a skew field fixing the column $\tgauche{(1, \, \hdots, \, 1)}$, we give formulas for a row vector fixed by this matrix.  The same techniques are applied to give noncommutative extensions of probabilistic properties of codes.
\end{abstract}
\maketitle

\section{Introduction} By a \textit{noncommutative stochastic matrix} we mean a square matrix $S = (a_{ij})_{1 \leq i, \, j \leq n}$ over a skew field, satisfying $\sum_{j = 1}^n \, a_{ij} = 1$ for any $i$; in other words, the row-sums are all equal to 1.  Equivalently the vector $\tgauche{(1, \, \hdots, \, 1)}$ is fixed by $S$.  We are answering here the following problem: find a row vector fixed by $S$.

\b

In the commutative case, formulas are known.  They occur in probability theory, where this problem is relevant.  Indeed, it amounts to finding the stationary distribution of the finite Markov chain whose transition matrix is $S$.  See Appendix 1 for details.

\b

But this problem may also be considered as a problem of noncommmutative linear algebra: given a square matrix over some skew field, which has a given column as eigenvector for some eigenvalue, find a corresponding row vector.  It is easy to reduce this general problem to the previous one, where the eigenvector is $\tgauche{(1, \, \hdots, \, 1)}$ and the eigenvalue is 1.

\b

In order to give formulas, which necessarily involve inverses of elements of the skew field and thus may be undefined, we take a \textit{generic noncommutative stochastic matrix}: this is the matrix $(a_{ij})$ of noncommuting variables $a_{ij}$ subject only to the condition that this matrix fixes $\tgauche{(1, \, \hdots, \, 1)}$.

\b

We seek now a row vector fixed by the matrix.  We work in the free field generated by these variables (in the sense of Paul Cohn), which we call the \textit{stochastic free field.}  The formula giving the row eigenvector uses the theory of variable-length codes.  Considering the complete digraph on the set $\{1, \, \hdots, \, n\}$, let $M_i$ be the set of paths from $i$ to $i$.  This is a free monoid and its basis $C_i$ is a prefix code.  Let $P_i$ be the set of proper prefixes of $C_i$, that is, the paths starting from $i$ and not passing through $i$ again.  We identify $P_i$ with the noncommutative power series which is equal to the sum of all the words in $P_i$ and we still denote this series by $P_i$. Then we show that the elements $P_i^{-1}$ can be evaluated in the stochastic free field and that the vector $(P_1^{-1}, \, \hdots, \, P_n^{-1})$ is fixed by our matrix; moreover, the $P_i^{-1}$ sum to 1, hence they form a kind of noncommutative limiting probability.  See Theorem \ref{thm: P_i sont definies dans un corps libre} and Example 1 to have a flavor of the result.

\b

The second part of the article deals with general variable-length codes, not necessarily prefix. One motivation is the fact that the proofs are quite similar.  The other motivation is that we obtain noncommutative generalization of well-known probabilistic results in the theory of codes, mostly due to Schützenberger (see \cite{BP86} and \cite{BPR}), who generalized the recurrent events of Feller.

\b

In Appendix 1, we review the commutative case.  In Appendix 2, we show how the theory of quasideterminants may be used to obtain our results on noncommutative stochastic matrices.

\b

\textit{Acknowledgments}

\n Thanks are due to Persi Diaconis and George Bergman for useful references; the article was improved by many suggestions of the latter.

\section{Basics}

\subsection{Langages and codes}
A \textit{language} is a subset of a free monoid $A^*$; the latter is generated by the \textit{alphabet} $A$.  A language is \textit{rational} if it is obtained from finite languages by the operations (called rational) \textit{union}, \textit{product} (concatenation) and \textit{star}.  The product of two languages $L_1L_2$ is $\{w_1w_2 \, \mid \, w_1 \in L_1, \, w_2 \in L_2\}$, and the star of $L$ is $L^* = \{w_1 \hdots w_n \, \mid \, w_i \in L, \, n \geq 0\} = \bigcup_{n \geq 0} \, L^n$.

\b

It is well-known that rational languages may be obtained by using only \textit{unambiguous rational operations}; these are: disjoint union, unambiguous product (meaning that if $w \in L_1 L_2$, then $w$ has a unique factorization $w = w_1 w_2$, $w_i \in L_i )$ and the star $L^*$ restricted to languages which are \textit{codes},  or equivalently bases of free submonoids of $A^*$.

\subsection{Series}

By a \textit{series} we mean an element of the $\mathbb{Q}$-algebra of noncommutative series $\mathbb{Q} \langle\langle A \rangle\rangle$, where $A$ is a set of noncommuting variables.  A \textit{rational} series is an element of the least subalgebra of $\mathbb{Q} \langle\langle A \rangle\rangle$, which contains the $\mathbb{Q}$-algebra of noncommutative polynomials $\mathbb{Q} \langle A \rangle$, and which is closed under the operation
$$S \mapsto S^* = \sum_{n = 0}^{\infty} \, S^n = (1 - S)^{-1},$$
which is defined if $S$ has zero constant term.  We denote by $\mathbb{Q} \langle\langle A \rangle\rangle^{rat}$ the $\mathbb{Q}$-algebra of rational series.  Each such series has a \textit{$*$- rational expression}: this is a well-formed expression involving scalars, letters (elements of $A$), products and star operations, the latter restricted to series with 0 constant term.  We say that a $*$-rational expression is \textit{positive} if the scalars involved are all $\geq 0$.

\b

Let $L$ be a rational language. Since $L$ may be obtained by unambiguous rational expressions, it follows that its \textit{characteristic series} $\sum_{w \in L} \, w \in \mathbb{Q} \langle\langle A \rangle\rangle$ is rational.  We shall identify a language and its characteristic series.  For all this, see \cite{BC88} or \cite{BC88bis}.

\subsection{Free fields}

The ring $\mathbb{Q} \langle\langle A \rangle\rangle^{rat}$ contains $\mathbb{Q} \langle A \rangle$; it is not a skew field.  There exist skew fields containing $\mathbb{Q} \langle A \rangle$.  Among them is the so called \textit{free field}.  We denote it $\mathcal{F}$.  It is generated by $\mathbb{Q} \langle A \rangle$ and has the following universal property (which characterizes it): for each $\mathbb{Q}$-algebra homomorphism $\mu: \mathbb{Q}\langle A \rangle \rightarrow D$, where $D$ is a skew field, there exists a $\mathbb{Q}$-subalgebra $O_{\mu}$ of $\mathcal{F}$ and a homomorphism $\overline{\mu}: O_{\mu} \rightarrow D$, extending $\mu$ and such that
$$f \in O_{\mu}, \q \overline{\mu} f \neq 0 \R f^{-1} \in O_{\mu}.$$

\b

The free field $\mathcal{F}$ is also characterized by the following property: say that a square matrix $M \in \mathbb{Q} \langle A \rangle^{n \times n}$ is \textit{full} if there exists no factorization $M = PQ$, $P \in \mathbb{Q} \langle A \rangle^{n \times r}$, $Q \in \mathbb{Q} \langle A \rangle^{r \times n}$, with $r < n$.  Then the square matrices over $\mathbb{Q} \langle A \rangle$ which are invertible in $\mathcal{F}$ are exactly the full matrices.  See \cite{Cohn71}.

\b

By a \textit{rational expression} over $\mathbb{Q} \langle A \rangle$ we mean a well-formed expression involving elements of $\mathbb{Q} \langle A \rangle$ and the operations sum, product and inversion.  Such an expression can be naturally evaluated in the free field $\mathcal{F}$, provided one never inverts 0.  For example, $(a + b^{-1})^{-1}$ can be evaluated in the free field whereas $(ab- (b^{-1}a^{-1})^{-1})^{-1}$ cannot. If rational expression can be evaluated in the free field, we say it is \textit{evaluable}.

\b

We shall consider also rational expressions over any skew field $D$, and say that such an expression is \textit{evaluable} if it can be evaluated without inversion of 0. If the elements of $D$ appearing in the rational expression are actually in a subring $R$ of $D$, we say that the expression is \textit{over $R$}.

\b

There is a canonical embedding of $\mathbb{Q} \langle\langle A \rangle\rangle^{rat}$ into $\mathcal{F}$, which can be seen as follows: let $S$ be any rational series; it has a $*$-rational expression; replace in it the operation $T^*$ by $(1 - T)^{-1}$; then one obtains a rational expression in $\mathcal{F}$, which is evaluable and represents the image of $S$ under the embedding $\mathbb{Q} \langle\langle A \rangle\rangle^{rat} \hookrightarrow \mathcal{F}$.  Thus, each rational language and each rational series is an element of the free field.  See \cite{Fliess}. In this way, each $*$-rational expression is equivalent to a evaluable rational expression over $\mathbb{Q} \langle A \rangle$.

\b

In the sequel, we use the notation $x^*$ for $(1-x)^{-1}$, when $x$ is in a ring and $1-x$ is invertible.

\subsection{The derivation $\lambda$ of the free field}

There is a unique derivation $\lambda$ of $\mathbb{Q} \langle A \rangle$ such that $\lambda(a)=a$, for any $a \in A$.  It maps each word $w \in A^*$ onto $|w|w$, where $|w|$ is the length of $w$. It has a unique extension to the free field $\mathcal{F}$, which we still denote $\lambda$.  Indeed, this follows from $Th. 7.5.17, p. 451$ in \cite{Cohn05}; see also \cite{Cohn00}.

\subsection{$D[t]$ and other rings}

Let $D$ be a skew field and $t$ be a central variable.  It is well-known that the ring of polynomials in $t$ over $D$ is a left and right Euclidean ring, and thus an Ore ring.  It has a field of fractions $D(t)$, each element of which is of the form $PQ^{-1}$ and $R^{-1}S$ for suitable $P, Q, R, S$ in $D[t]$.  The ring of series in $t$ over $D$ is denoted $D[[t]]$.  It is contained in the skew field of Laurent series $D((t))$. The latter also contains canonically $D(t)$, and we may identify $D(t)$ with a subset of $D((t))$.  A series $S \in D[[t]]$ is called \textit{rational} if $S$ is in $D(t)$.  The ring of rational series is denoted by $D[[t]]^{rat}$.  Thus $D[[t]]^{rat} = D[[t]] \cap D(t)$.

\b

Each polynomial $P\in D[t]$ may uniquely be written $P = (1-t)^n Q$, with $n \in \mathbb{N}$, $Q \in D[t]$ and $Q(1) \neq 0$.  Thus if $S \in D(t)$, one has $S = (1-t)^nQR^{-1}$, with $n \in  \mathbb{Z}$, $Q, \, R \in D[t]$ and $Q(1), \, R(1) \neq 0$.  We say that $S$ is \textit{evaluable at $t = 1$} if $n \geq 0$, and in this case, its \textit{value at $t = 1$} is $Q(1) R(1)^{-1}$ if $n = 0$, and 0 if $n \geq 1$.  This value is a well-defined element of $D$, which does not depend on the fraction chosen to represent $S$.

\b

We extend this to matrices: a matrix over $D(t)$ is said to be \textit{evaluable at $t = 1$} if all his entries are, and then its value at $t = 1$ is defined correspondingly.

\b

Consider a rational expression $E(t)$ over the skew field $D(t)$. We obtain a rational expression over the skew field $D$ by putting $t = 1$ in $E(t)$.  Suppose that the rational expression $E(1)$ obtained in this way is evaluable in $D$ and evaluates to $\alpha \in D$; then the rational expression $E(t)$ is evaluable in $D(t)$, evaluates to an element $P(t)Q(t)^{-1}$ in $D(t)$, with $P, \, Q \in D[t]$, and $PQ^{-1}$ is evaluable at $t = 1$ with value $\alpha \in D$.  The standard details are left to the reader.

\subsection{Central eigenvalues of matrices over a skew field}

Let $M$ be a square matrix over $D$.  Then $1 - Mt$ is invertible over $D[[t]]$, hence over $D((t))$.  Since $D(t)$ is a skew field, contained in $D((t))$, and containing the coefficients of $1 - Mt$, the coefficients of its inverse $(tM)^* = (1-tM)^{-1}$ lie also in $D(t)$ and finally in $D[[t]]^{rat}$.  Recall that a square matrix over a skew field is left singular (that is, has a nontrivial kernel when acting at the left on column vectors) is and only if it is right singular.  Thus $M$ has an eigenvector for the eigenvalue 1 at the left if and only if it holds on the right.

\b

By the \textit{multiplicity} of the eigenvalue 1 of $M$ we mean the maximum of the nullity (that is, dimension of kernel) of the positive powers of $1 - M$.  Observe that this coincides with the usual multiplicity (of 1 in the characteristic polynomial) if $D$ is commutative.  Note that the same properties hold for any nonzero central eigenvalue $\lambda$ by considering $1- \lambda^{-1}M$; we treat only the case $\lambda = 1$ for the future application.

\begin{lem} \label{lemme concernant la valeur propre 1 de I-tm}
Let $M$ be a square matrix over the skew-field $D$ and $t$ be a central variable.
\newcounter{marqueur11}
\begin{list}{(\roman{marqueur11})}{\usecounter{marqueur11}}
\item $M$ has the eigenvalue $1$ if and only if $(1 - tM)^{-1}$ is not evaluable at $t = 1$. \\
\item If $M$ has the eigenvalue $1$ with multiplicity $1$, then $(1-t)(1-tM)^{-1}$ is evaluable at $t = 1$, is nonnull and its rows span the eigenspace for the eigenvalue $1$.
\end{list}
\end{lem}

\n \textit{Proof}
\b
\n (i) Suppose that $M$ has the eigenvalue 1.  Then $M$ is conjugate over $D$ to a matrix of the form $N = \begin{bmatrix}
1 & 0 \\
P & Q
\end{bmatrix}$, where $Q$ is square.  Then, computing in $D[[t]]$, we have $(1-tN)^{-1} = \begin{bmatrix}
    (1-t)^{-1} & 0 \\
    \times & \times
    \end{bmatrix}.$
    This is clearly not evaluable for $t = 1$, and therefore $(1 - tM)^{-1}$ is also not evaluable for $t = 1$.

    Conversely, suppose $(1-tM)^{-1}$ is not evaluable for $t = 1$.  Then, we have $(1-tM)^{-1} = ((1-t)^{n_{ij}} \, P_{ij} / Q_{ij})_{i, \, j}$, with $P_{ij}, \, Q_{ij} \in D[t]$, $P_{ij}(1), \, Q_{ij}(1) \neq 0$, $n_{ij} \in \mathbb{Z}$ and some $n_{ij} < 0$.  Let $-n$ be the minimum of the $n_{ij}$.  Then $n > 0$ and $(1-t)^n \, (1-tM)^{-1}$ is evaluable at $t = 1$ and its value $P$ at $t = 1$ is nonnull.  Now, we have
    $$(1-tM)^{-1} = 1 + (1-tM)^{-1}tM,$$
    thus
    $$(1-t)^n (1-tM)^{-1} = (1-t)^n + (1-t)^n(1-tM)^{-1}tM.$$
    Since $n > 0$, we obtain for $t = 1$:
    $$P = PM,$$
    which shows that $M$ has the eigenvalue 1, since each row of $P$ is fixed by $M$ and $P \neq 0$.

  \n $(ii)$ We write as before $N = \begin{bmatrix}
1 & 0 \\
P & Q
\end{bmatrix}$, where $N$ is conjugate to $M$ over $D$.  Then
$$(1-tN)^{-1} = (tN)^* = \begin{bmatrix}
t^* & 0 \\
(tQ)^*tPt^* & (tQ)^*
\end{bmatrix}.$$
We claim that $(tQ)^*$ is evaluable at $t = 1$. Indeed, otherwise, by $(i)$, $Q$ has the eigenvalue 1 and is conjugate to a matrix $N = \begin{bmatrix} 1 & 0 \\ R & S \end{bmatrix}$, $S$ square.  Then $M$ is conjugate to
$N = \begin{bmatrix}
1 & 0 & 0 \\
\times & 1 & 0 \\
\times & R & S
\end{bmatrix}$
and the square of $1-M$ has nullity $\geq 2$, contradiction.

Now, we see that

$$(1-t)(tN)^* = \begin{bmatrix}
1 & 0 \\
(tQ)^*tP & (1-t)(tQ)^*
\end{bmatrix}$$
is evaluable at $t = 1$
and that its value at $t=1$ is nonnull. Thus, by the first part of the proof, its rows span the eigenspace for the eigenvalue 1.
\begin{flushright}
$\square$
\end{flushright}

\subsection{Rational series in one variable}

Let $R$ be a ring and $t$ a central variable.  In the ring of formal power series $R[[t]]$, we consider the subring $R[[t]]^{rat}$, which is the smallest subring containing $R[t]$ and closed under inversion.  If $R$ is a skew field, then $R[[t]]^{rat}$ canonically embeds into the skew field $R(t)$.  If $R \rightarrow S$ is a ring homomorphism, then it induces a ring homomorphism $R[[t]]^{rat} \rightarrow S[[t]]^{rat}$ fixing $t$.

\section{Generic noncommutative stochastic matrices}

\subsection{Generic matrices}

Let $M = (a_{ij})_{1 \leq i, \, j \leq n}$ be a \textit{generic noncommutative matrix}; in other words, the $a_{ij}$ are noncommuting variables.  We denote by $\mathcal{F}$ the corresponding free field.  Associated to $M$ is the matrix $S$: it is the same matrix, but this time we assume that the $a_{ij}$ are noncommuting variables subject to the \textit{stochastic identities}
\begin{equation}
\forall i = 1, \, \hdots, \, n, \, \sum_{j = 1}^n \, a_{ij} = 1. \label{relations stochastiques}
\end{equation}
In other words, the row sums of $S$ are equal to 1; equivalently, $S$ has $\tgauche{(1, \, \hdots, \, 1)}$ as column eigenvector with the eigenvalue 1.  We call $S$ a \textit{generic noncommutative stochastic matrix}.  The algebra over $\mathbb{Q}$ generated by its coefficients (hence subject to the relations (\ref{relations stochastiques})) is a free associative algebra, since it is isomorphic with $\mathbb{Q} \langle a_{ij}, \, i \neq j\rangle$.  Indeed, we may eliminate $a_{ii}$ using (\ref{relations stochastiques}).  We denote this algebra by $\mathbb{Q} \langle a_{ij} / (\ref{relations stochastiques}) \rangle$, referring to the relations (\ref{relations stochastiques}). Hence there is a corresponding free field, which we call the \textit{stochastic free field}, denoted $\mathcal{S}$.

\subsection{Existence of elements and identities in the stochastic free field}

We want to verify that certain rational expressions make sense in the stochastic free field $\mathcal{S}$.  For example, anticipating on the example to come, we want to show that $(1 + bd^*)^{-1} = (1 + b(1-d)^{-1})^{-1}$ makes sense in $\mathcal{S}$ (hence under the hypothesis $a + b = c+ d = 1)$.  It is necessary to take care of this existence problem, since otherwise, one could invert 0 (and our proved identities will be meaningless). The idea is to prove the existence of certain specializations of the variables, compatible with the identities in $\mathcal{S}$ (identities (\ref{relations stochastiques}) above), such that the specialized rational expression makes sense.  In our example, we could take $b = 0$:  then $bd^*$ specializes to 0 and $1 + bd^*$ to 1, hence $(1 + bd^*)^{-1}$ is evaluable under the specialization.  A fortiori, since $\mathcal{S}$ is a free field, $(1 + bd^*)^{-1}$ is evaluable in $\mathcal{S}$.

\b

By a \textit{Bernouilli morphism} we mean a $\mathbb{Q}$-algebra morphism $\pi$ of the free associative algebra $\mathbb{Q}\langle a_{ij} \rangle$ into $\mathbb{R}$ such that
\newcounter{marqueur1}
\begin{list}{(\roman{marqueur1})}{\usecounter{marqueur1}}
\item for any $i = 1, \, \hdots, \, n$, $\sum_{j = 1}^n \, \pi (a_{ij}) = 1$;
\item $\pi (a_{ij}) > 0$, for any $i, \, j = 1, \, \hdots, \, n$.
\end{list}
Clearly, such a morphism induces naturally a $\mathbb{Q}$-algebra morphism form $\mathbb{Q} \langle a_{ij} / (\ref{relations stochastiques}) \rangle$ into $\mathbb{R}$.

\begin{lem}\label{lemme: existence d'un sous-corps}There exists a subring $\mathcal{S}_{\pi}$ of the stochastic free field $\mathcal{S}$ such that
\newcounter{marqueur2}
\begin{list}{(\roman{marqueur2})}{\usecounter{marqueur2}}
\item $\mathcal{S}_{\pi}$ contains $\mathbb{Q} \langle a_{ij} / (\ref{relations stochastiques}) \rangle$;
\item there is an extension of $\pi$ to $\mathcal{S}_{\pi}$ (we still denote it by $\pi$);
\item if $f \in \mathcal{S}_{\pi}$ and $\pi(f) \neq 0$, then $f^{-1} \in \mathcal{S}_{\pi}$.
\end{list}
\end{lem}

\n \textit{Proof}
This is a consequence of the fact that $\mathcal{S}$ is a free field, corresponding to the free associative algebra $\mathbb{Q} \langle a_{ij} \, /(\ref{relations stochastiques}) \rangle$, hence is the universal field of fractions of $\mathbb{Q} \langle a_{ij} \, /(\ref{relations stochastiques}) \rangle$.  This implies that there exists a specialization $\mathcal{S} \rightarrow \mathbb{R}$ extending $\pi: \mathbb{Q} \langle a_{ij} / (\ref{relations stochastiques}) \rangle \rightarrow \mathbb{R}$, and the lemma follows from \cite{Cohn71} 7.2 and Cor. 7.5.11.
\begin{flushright}
$\square$
\end{flushright}

\begin{cor} \label{corollaire: mesure de Bernouilli => expression rationnelle existe dans S}
Suppose that $\pi$ is a Bernouilli morphism and that $S = \sum_{w \in L} \, w$, where $L$ is a rational subset  of the free monoid $\{a_{ij}\}^*$ such that $\sum_{w \in L} \, \pi(w) < \infty$.  Then any positive $*$-rational expression for $S$ is evaluable in the stochastic free field $\mathcal{S}$.
\end{cor}

\b

\n \textit{Proof}
This is proved inductively on the size of the rational expression for $S$.  Note that for each subexpression and corresponding series $S'$, $\pi(S')$ converges and is $> 0$.  Hence, we apply inductively the lemma and see that for each subexpression, the corresponding element is in $\mathcal{S}_{\pi}$.
\begin{flushright}
$\square$
\end{flushright}

\begin{lem} \label{lemme: S definie dans varfi, S definie dans F}Let $S$ be a $*$-rational series in $\mathbb{Q} \langle\langle a_{ij}\rangle\rangle$ having a $*$-rational expression which is evaluable in $\mathcal{S}$.  Then it is evaluable in the free field $\mathcal{F}$.  If moreover $S = 0$ in $\mathcal{F}$, then $S = 0$ in $\mathcal{S}$.
\end{lem}

\n \textit{Proof}
There exists a specialization $\mathcal{F} \rightarrow \mathcal{S}$, since $\mathcal{F}$ is the universal field of fractions of $\mathbb{Q} \langle a_{ij}\rangle$, see \cite{Cohn71} chapter 7.  Hence there is a subring $H$ of $\mathcal{F}$ and a surjective $\mathbb{Q}$-algebra morphism $\sigma: H \rightarrow \mathcal{S}$ such that: $\forall \, f \in H$, $\sigma{f} \neq 0 \R f^{-1} \in H$, and such that $H$ contains $\mathbb{Q} \langle a_{ij}\rangle$.

\b

We may therefore prove, by induction on the size of the rational expression, that $S$ exists in $H$ and that $\sigma(S)$ is the element of $\mathcal{S}$ defined by the rational expression.  It follows that, if $S = 0$ in $\mathcal{F}$, then $S = 0$ in $\mathcal{S}$.
\begin{flushright}
$\square$
\end{flushright}

\subsection{Paths}

Each path in the complete directed graph with set of vertices $\{1, \, \hdots, \, n\}$ defines naturally an element of the free associative algebra $\mathbb{Q} \langle a_{ij}\rangle$, hence of the free field $\mathcal{F}$.  This is true also for each rational series in $\mathbb{Q} \langle \langle a_{ij} \rangle \rangle$.

\b

We define now several such series.  First, consider the set of nonempty paths $i \rightarrow i$ which do not pass through $i$; we denote by $C_i$ the sum in $\mathbb{Q} \langle \langle a_{ij} \rangle \rangle$ of all the corresponding words.  It is classically a rational series, and thus defines an element of the free field $\mathcal{F}$. Now, let $P_i$ be the sum of the paths (that is, the corresponding words) from $i$ to any vertex $j$, which do not pass again through $i$; this set of words is the set of proper prefixes of the words appearing in $C_i$.  Likewise, $P_i$ defines an element of $\mathcal{F}$.

\b

\n \textbf{Example 1.} \label{exemple 1}
$M = \begin{bmatrix}
a & b \\
c & d
\end{bmatrix}$.  The graph is
\begin{figure}[H]
  \begin{center}
    \unitlength=3pt
    \begin{picture}(25, 5)
    \gasset{Nw=5,Nh=5,Nmr=2.5,curvedepth=4}
    \thinlines
    \node(A1)(0,0){$1$}
    \node(A2)(25,0){$2$}
    \drawloop[loopangle=180](A1){$a$}
    \drawloop[loopangle=0](A2){$d$}
    \drawedge(A1,A2){$b$}
    \drawedge(A2,A1){$c$}
    \end{picture}
  \end{center}
\end{figure}
Then
\begin{eqnarray*}
C_1 & =  a + b d^* c, \qquad C_2 & =  d + c a^* b, \\
P_1 & =  1 + b d^*, \qquad P_2 & =  1 + c a^*.
\end{eqnarray*}

\subsection{Results}

\begin{thm} \label{thm: P_i sont definies dans un corps libre}
The elements $P_i$ can be evaluated in the stochastic free field $\mathcal{S}$ and $(P_1^{-1}, \, \hdots, \, P_n^{-1})$ is a left eigenvector of the noncommutative generic stochastic matrix $S$.  Moreover, in $\mathcal{S}$,
\newcounter{marqueur3}
\begin{list}{(\roman{marqueur3})}{\usecounter{marqueur3}}
\item $\displaystyle{\sum_{i = 1}^n} \, P_i^{-1} = 1;$
\item $C_i$ can be evaluated in $\mathcal{S}$ and is equal to $1$;
\item $\lambda(C_i)$ can be evaluated in $\mathcal{S}$ and is equal to $P_i$.
\end{list}
\end{thm}

\b

Here $\lambda$ is the unique derivation of the free field $\mathcal{F}$ which extends the identity on the set $\{a_{ij}\}$.

\n \textbf{Exemple 1.} (continued)  We verify that $(P_1^{-1}, \, P_2^{-1}) \, \begin{bmatrix}
a & b \\
c & d
\end{bmatrix} = (P_1^{-1}, \, P_2^{-1})$.
It is enough to show that
$P_1^{-1} \, a + P_2^{-1} \, c = P_1^{-1}$.  This is equivalent to
\begin{eqnarray*}
P_2^{-1} \, c & = & P_1^{-1} (1 - a) \\
\LR c^{-1} \,P_2 & = & a^* P_1 \\
\LR c^{-1} + a^* & = & a^* + a^* b d^*.
\end{eqnarray*}
Now, we take the stochastic identities:
\begin{eqnarray*}
a + b & = & 1 \R 1 - a = b \R a^* = b^{-1} \R a^*b = 1, \\
c + d & = & 1 \R d^* = c^{-1}.
\end{eqnarray*}
Thus, we may conclude.
\newcounter{marqueur4}
\begin{list}{(\roman{marqueur4})}{\usecounter{marqueur4}}
\item Similarly:
\begin{eqnarray*}
P_1^{-1} + P_2^{-1} & = & 1 \\
\LR P_2 + P_1 & = & P_1 P_2 \\
\LR 2 + bd^* + ca^* & = & 1 + bd^* + ca^* + bd^*ca^*
\end{eqnarray*}
and we conclude since $d^*c = ba^* = 1$.
\item In $\mathcal{S}$, $C_1 = a + bd^*c = a + b = 1.$
\item In $\mathcal{F}$, $\lambda(C_1) = a + bd^*c + b \lambda(d^*)c + bd^*c,$ since $\lambda$ is a derivation such that $\lambda(b) = b$ and $\lambda(c) = c$.  Now
$$\lambda(d^*) = \lambda((1-d))^{-1} = -(1-d)^{-1}\lambda(1-d) (1-d)^{-1} = d^*d d^*.$$
Thus, this time in $\mathcal{S}$,
\begin{eqnarray*}
\lambda(C_1) & = & a + 2 bd^*c + bd^*dd^*c \\
& = & a + 2b  +bd^*d \\
& = & a+b + b(1+d^*d) \\
& = & 1 + bd^*  = P_1.
\end{eqnarray*}
\end{list}

\subsection{Proof of the theorem}

\begin{lem} \label{Lemme: (1-t)(tS)^* well-defined}
Consider the matrix $(tS)^*$ in $\mathcal{S}(t)$.  Then $(1-t)(tS)^*$ is can be evaluated for $t=1$ and is nonzero.
\end{lem}

\n \textit{Proof}
By Lemma \ref{lemme concernant la valeur propre 1 de I-tm}, it is enough to show that $S$ has the eigenvalue 1 with multiplicity 1.  Now, by a change of basis over $\mathbb{Q}$ (replace the canonical basis of column vectors $e_1, \, \hdots, \, e_n$ by $e_1, \, \hdots, \, e_{n-1}, \, e_1 + \hdots + e_n$), we bring $S$ to the form
$$T = \begin{bmatrix}
N & 0 \\
\lambda & 1
\end{bmatrix},$$
where $n_{ij} = a_{ij} - a_{nj}$ for $1 \leq i, \, j \leq n - 1$ and $\lambda_j = a_{nj}$ for $j = 1, \, \hdots, \, n-1$.  We claim that $N-1$ is inversible in $\mathcal{S}$.  It is enough to show that it is full in $\mathbb{Q} \langle a_{ij} \, /(\ref{relations stochastiques}) \rangle$. Suppose that $N-1$ is not full: then $N-1 = PQ$, with $P, \, Q$ over $\mathbb{Q} \langle a_{ij} \, /(\ref{relations stochastiques}) \rangle$ of size $n \times (n-1)$ and $(n-1) \times n$.  By replacing $a_{nj}$ by 0 and $a_{ii}$ by $a_{ii}+1$, we find that the matrix $(a_{ij})_{1 \leq i, \, j \leq n-1}$ is nonfull over $\mathbb{Q} \langle a_{ij}, \, 1 \leq i, \, j \leq n-1 \rangle$, which is absurd, since it is a generic matrix.  Thus $N-1$ is inversible, and no power of it has a kernel.  Consequently, the positive powers of $T-1$ have all rank $n-1$.  Therefore the multiplicity of 1 as eigenvalue of $T$, hence of $S$, is 1.
\begin{flushright}
$\square$
\end{flushright}

\n \textit{Proof of Theorem \ref{thm: P_i sont definies dans un corps libre}}

Let us identify paths in the complete directed graph on $\{1, \, \hdots, \, n\}$, and corresponding words in the free monoid $\{a_{ij}\}^*$.  We identify also an infinite sum of paths with the corresponding series in $\mathbb{Q}\langle\langle a_{ij}\rangle\rangle$.  Let $P_{ij}$ denote the set of paths from $i$ to $j$ that do no pass through $i$ again.  We therefore have $P_i = \sum_j \, P_{ij}$.  Denote by $D(u_1, \, \hdots, \, u_n)$ the diagonal matrix whose diagonal elements are $u_1, \, \hdots, \, u_n$.  Observe that each path from $i$ to $j$ may be decomposed as the concatenation of a path from $i$ to $i$ (thus, an element of $C_i^*$) and a path from $i$ to $j$ that does not pass again through $i$ (thus, an element of $P_{ij}$).  Since $(M^*)_{ij}$ is the sum of all paths from $i$ to $j$, we obtain the identity in $\mathbb{Q} \langle\langle a_{ij }\rangle\rangle:(M^*)_{ij} = C_i^* \, P_{ij}$.  Thus we have the matrix identity: $M^* = D(C_1^*, \, \hdots, \, C_n^*) (P_{ij})$.  Now $P_{ii} = 1$ and $P_{ij}$ has no constant term.  Hence $(P_{ij})$ is invertible over $\mathbb{Q} \langle\langle a_{ij} \rangle\rangle$.

\b

Inverting, we obtain $D(C_1 - 1, \, \hdots, \, C_n - 1) = (P_{ij})(M-1)$, since $M^* = (1 - M)^{-1}$ and similarly $C_i^* = (1- C_i)^{-1}$.  If we multiply by the column vector $\gamma = \tgauche{(1, \, \hdots, \, 1)}$, we obtain $\tgauche{(C_1 - 1, \, \hdots, \, C_n - 1)} = (P_{ij}) \, (M-1) \, \gamma$.

\b

This equality holds in $\mathbb{Q} \langle\langle a_{ij} \rangle\rangle$, and  actually, in its subalgebra of rational series, since $C_i, \, P_{ij}$ are rational series.  Hence it holds in the free field $\mathcal{F}$.

\b

We also obtain, applying the derivation $\lambda$ of $\mathcal{F}$:
$$\tgauche{(\lambda(C_1), \, \hdots, \, \lambda(C_n))} = (\lambda(P_{ij})) \, (M - 1)  \, \gamma + (P_{ij}) \, M  \, \gamma.$$
Now, we claim that $C_i, \, P_{ij}$ and $\lambda(C_i)$ can be evaluated in the stochastic free field $\mathcal{S}$.  Thus, since $M \gamma = \gamma$ in $\mathcal{S}$, we obtain that in $\mathcal{S}:$
$$\tgauche{(C_1 - 1, \, \hdots, \, C_n - 1)} = 0, \q$$ and $$\q \tgauche{(\lambda(C_1), \, \hdots, \, \lambda(C_n))} = (P_{ij}) \, \gamma = \tgauche{(P_1, \, \hdots, \, P_n)},$$
\n which proves parts $(ii)$ and $(iii)$ of the theorem.

\b

In order to prove the claim, we take a Bernouilli morphism $\pi$.  Let $i$ be some element of $\{1, \, \hdots, \, n\}$ and consider the set $E$ of paths not passing through $i$.  Then $\pi(E) < \infty$ since the matrix $N$, which is obtained from $M$ by removing row and column $i$, satisfies: $\pi(N)$ has row sums $< 1$.  It follows that $\pi(C_i), \, \pi(P_{ij})$ are finite.  For $\lambda(C_i)$, it is easy to see inductively on the size of a rational expression of $C_i$ that, since $C_i$ can be evaluated in $\mathcal{S}$, so is $\lambda(C_i)$; one has simply to use the identity $\lambda(H^*) = H^* \, \lambda(H) \, H^*$.  Note also that $\pi(P_i) > 0$, hence $P_i$ is nonzero in $\mathcal{S}$, and $P_i^{-1}$ is an element in $\mathcal{S}$, by Corollary \ref{corollaire: mesure de Bernouilli => expression rationnelle existe dans S}.

\b

We now prove $(i)$. Let $Q_i$ denote the set of paths from 1 to some vertex, that do not pass by $i$; in particular, $Q_1 = 0$. Then, for any $i, \, j$, we have
$$(M^*)_{1i} \, P_i + Q_i = (M^*)_{1j} \, P_j + Q_j,$$
since both sides represent all the paths departing from 1.  Let $t$ be a central variable.  Replacing each path $w$ by $t^{|w|} \, w$ and writing correspondingly $P_1(t), \, \hdots, \, P_n(t)$, we obtain
$$(tM)_{1i}^* \, P_i(t) + Q_i(t) = (tM)_{1j}^* \, P_j(t) + Q_j(t).$$
This holds in $\mathbb{Q}\langle A \rangle[[t]]$ and actually in its subalgebra of rational elements $\mathbb{Q}\langle A \rangle[[t]]^{rat}$.  Now, we have canonical homomorphisms (see 2.5 and 2.7)
$$\mathbb{Q}\langle A \rangle[[t]]^{rat} \rightarrow \mathbb{Q}\langle A / (\ref{relations stochastiques}) \rangle[[t]]^{rat} \rightarrow \mathcal{S}[[t]]^{rat} \rightarrow \mathcal{S}(t).$$

\b

The composition maps the matrix $M$ onto $S$.  Hence, we have in $\mathcal{S}(t)$
$$(tS)_{1i}^* \, P_i(t) + Q_i(t) = (tS)_{1j}^* \, P_j(t) + Q_j(t),$$
where we keep the notation $P_i(t) \in  \mathcal{S}(t)$ for the image under the composition.  Observe that $P_i$, by Cor. \ref{corollaire: mesure de Bernouilli => expression rationnelle existe dans S}, has a rational expression which can be evaluated in $\mathcal{S}$.  Hence $P_i(t)$ can be evaluated for $t = 1$ and equal to $P_i$.  Similarly, $Q_i(t)$ can be evaluated for $t = 1$ and evaluates to $Q_i$.

\b

Multiply the last equality by $1 - t$.  By Lemma \ref{Lemme: (1-t)(tS)^* well-defined}, $(1-t) \, (tS)_{1i}^*$ can be evaluated for $t = 1$ and is equal to $\alpha_i$ say. Thus, we obtain
$$\alpha_i \, P_i = \alpha_j \, P_j.$$
Now $(tS)^* = 1 + (tS)^* \, tS$, so that, putting $t = 1$, we obtain that each row of $(1-t)(tS)^* \, \b |_{t = 1}$ is fixed by $S$. In particular, so is $(\alpha_1, \, \hdots, \, \alpha_n)$.  Since by Lemma \ref{Lemme: (1-t)(tS)^* well-defined}, $(1-t)(tS)^* \, \b |_{t = 1}$ is nonzero, some row of it is nonzero, and by symmetry, each row is nonzero.  Hence, since we already know that each $P_i$ is nonzero in $\mathcal{S}$, we see that each $\alpha_i$ is $\neq 0$.  Thus, since $P_i^{-1} \alpha_i^{-1} = P_1^{-1} \alpha_1^{-1}$,
$$(P_1^{-1}, \, \hdots, \, P_n^{-1}) = P_1^{-1} \, \alpha_1^{-1} (\alpha_1, \, \hdots, \, \alpha_n),$$
which shows that $(P_1^{-1}, \, \hdots, \, P_n^{-1})$ is fixed by $S$.

\b

Now $\sum_{i = 1}^n \, (M^*)_{1i} = M_{11}^* \, P_1$, since both sides represent the paths departing from 1.  Thus we deduce that $\sum_{i = 1}^n \, \alpha_i = \alpha_1 \, P_1$
in $\mathcal{S}$, by the same technique as above.  Thus
$$\sum_{i = 1}^n \, P_i^{-1} = \sum_{i = 1}^n \, P_1^{-1} \, \alpha_1^{-1} \, \alpha_i = 1.$$
\begin{flushright}
$\square$
\end{flushright}

\section{Unambiguous automata}

\subsection{Unambiguous automata}

An \textit{unambiguous automata} is equivalent to a multiplicative homomorphism $\mu$ from the free monoid $A^*$ into $\mathbb{Q}^{n \times n}$ such that each matrix $\mu w$, $w \in A^*$, has entries in $\{0, \, 1\}$.  This may be expressed by associating to $\mu$ the directed graph with edges labelled in $A$ with vertices $1, \, \hdots, \, n$, and edges $\xymatrix @!0 @C=1cm {i \ar[r]^{a} & j}$ if and only if $(\mu a)_{ij} = 1$.  Then the non-ambiguity means that
for any vertices $i, \, j$ and any word $w$, there is at most one path from $i$ to $j$ labelled $w$ (the label of a path is the product of the label of the edges).  The \textit{matrix of the automaton} is by definition $M = \sum_{a \in A} \, a \, \mu a$.

\b

We say that the unambiguous automaton is \textit{complete} if the zero matrix does not belong to the monoid $\mu A^*$. Equivalently, for each word $w$ there is some path labelled $w$.  We say that the automaton is \textit{transitive} is the underlying graph is strongly connected. This means that for any vertices $i, \, j$, there is some path $i \rightarrow j$; equivalently, $(\mu w)_{ij} \neq 0$ for some word $w$.

\b

The monoid $\mu A^*$ is finite.  Hence it has a unique minimal ideal $I$.  There is a \textit{rank} function on $\mu A^*$, and the elements of minimum rank are precisely the elements of $I$.  Since $\mu A^* \subseteq \{0, \, 1\}^{n \times n}$, the rows of an element in $\mu A^*$ are ordered by inclusion (by identifying a subset of  $\{1, \, 2, \, \hdots, \, n\}$ and its characteristic row vector).  It is shown that the nonzero rows of elements of the minimal ideal are precisely the maximal rows of elements of $\mu A^*$.  Similarly for columns.  The ideal $I$ is the disjoint union of the minimal right (resp. left) ideals of $\mu A^*$, and the intersection of a minimal left and a minimal right ideal is a group. For this, see \cite{BP86} Chapter $VI$, and \cite{BPR} Chapter $VI$, especially Exercice 3.4 and also \cite{BCKP}.

We shall use the following result

\begin{prop}
Let $c$ be a maximal column and $R$ be the sum of the distinct rows of some element in the minimal ideal.  Then $Rc = 1$.
\end{prop}

\b

\n \textit{Proof}
There exist $x, \, y$ in $I$ such that $c$ is a column of $x$ and $R$ is the sum of the distinct rows of $y$.  The element $xy$ is in $I$ and belongs therefore to a group with neutral element $e$, say.  Then $e$ has a column-row decomposition $e=st$, where $s$ (resp. $t$) is a $n \times r$ (resp. $r \times n$) matrix with entries in $\{0, \, 1\}$, with $r$ the minimal rank of $\mu A^*$, where $ts = I_r$ (the identity matrix), and where the set of nonzero rows of $e$ is the set of rows of $t$, which has distinct rows, and similarly for the columns of $s$ (see \cite{BP86} Prop. IV.3.3 or \cite{BPR}, Prop. VI.2.3).

\b

Now, $xM$ is a minimal right ideal of $\mu A^*$, containing $e$, hence $xM = eM$ and therefore $x = em = stm$.  Hence $c$ is a sum of columns of $s$, and since $c$ is a maximal column, $c$ is a column of $s$.  Similarly, $y = nst$ and each nonzero row of $y$ is a row of $t$. We have also $e=n'y$, hence each nonzero row of $t$, being a row of $e$, is a row of $y$.  Thus $R$ is the sum of the rows of $t$: $R = \lambda t$, with $\lambda = (1, \, \hdots, \, 1)$.  Finally $Rc = \lambda t c$ and since $ts = I_r$ and $c$ is a column of $s$, $tc$ is a column of $I_r$ and $\lambda tc = 1.$
\begin{flushright}
$\square$
\end{flushright}

\b

\n \textbf{Example 2.}

The unambiguous automaton is
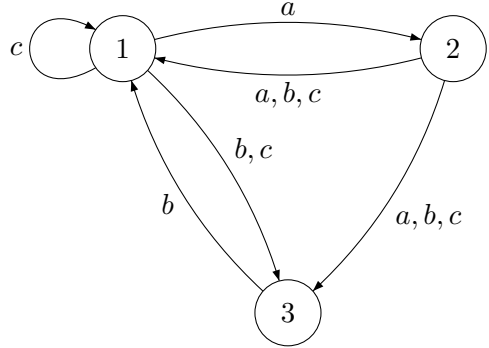
\begin{figure}[H]
  \begin{center}
    \unitlength=5pt
    \begin{picture}(25, 29)
    \gasset{Nw=5,Nh=5,Nmr=2.5,curvedepth=3}
    \thinlines
    \node(A1)(0,20){$1$}
    \node(A2)(25,20){$2$}
    \node(A3)(12.5,0){$3$}
    \drawloop[loopangle=180](A1){$c$}
    \gasset{curvedepth=2}
    \drawedge(A1,A2){$a$}
    \drawedge(A2,A1){$a, b, c$}
    \drawedge(A3,A1){$b$}
    \drawedge(A1,A3){$b, c$}
    \drawedge(A2,A3){$a, b, c$}
    \end{picture}
  \end{center}
\end{figure}
The associated representation $\mu$ is defined by
$$\mu a = \begin{bmatrix}
0 & 1 & 0 \\
1 & 0 & 1 \\
0 & 0 & 0
\end{bmatrix}, \q \mu b = \begin{bmatrix}
0 & 0 & 1 \\
1 & 0 & 1 \\
1 & 0 & 0
\end{bmatrix}, \q \mu c = \begin{bmatrix}
1 & 0 & 1 \\
1 & 0 & 1 \\
0 & 0 & 0
\end{bmatrix}.$$
The matrix of the automaton is
$$M = \begin{bmatrix}
c & a & b + c \\
a + b + c & 0 & a + b + c\\
b & 0 & 0
\end{bmatrix}$$
Idempotents in the minimal ideal are for example $\mu c$ and $\mu ba = \begin{bmatrix}
0 & 0 & 0 \\
0 & 1 & 0 \\
0 & 1 & 0
\end{bmatrix}$.
The maximal rows are $(1, \, 0, \, 1)$ and $(0, \, 1, \, 0)$ and the maximal columns are $\tgauche{(1, \, 1, \, 0)}$ and $\tgauche{(0, \, 1, \, 1)}$.

\subsection{Codes}

Recall that a code is the basis of some free submonoid of the free monoid.  Given an unambiguous automaton with associated representation $\mu$, and some vertex $i$, the language $\{ w \in A^* \, \mid \, (\mu w)_{ii} = 1\}$ is a free submonoid of $A^*$; we denote by $C_i$ its unique basis, which is therefore a code.  This code is moreover rational.  Explicitly, $C_i$ is the set of labels of paths $i \rightarrow i$ which do not contain $i$ as internal vertex.  Note that $C_i$ is a rational code and that each rational code is obtained in this way. We shall use also the set $P_i$ of labels of paths starting at $i$ and not passing again through $i$. See \cite{BP86}.

\n \textbf{Example 2.} (continued)
Write $A = a + b + c$, then
\begin{eqnarray*}
C_1 & = & c + aA(1 + Ab) + (b+c)b, \\
C_2 & = & A(c + b^2 + cb)^*a + Ab(c + b^2 + cb)^* a, \\
C_3 & = & b(c + aA)^* (b + c + aA), \\
P_1 & = & 1 + a + aA + b + c, \\
P_2 & = & 1 + A(c + b^2 + cb)^* (1 + b + c) + A(1+b (c + b^2  +cb)^* (b + c)) + Ab(c + b^2  +cb)^*, \\
P_3 & = & 1 + b(c + aA)^* (1 + a).
\end{eqnarray*}

We  shall use the following property of rational maximal codes: let $C$ be such a code; then there exists rational languages $P, \, S, \, F$, whose elements are factors of words of $C$, such that in $\mathbb{Q} \langle\langle A \rangle\rangle$
$$A^* = S C^*P + F.$$
Moreover $1 \not \in F$, $1 \in S$, $1 \in P$. This property is proved in \cite{BPR} Lemma XII.4.3 for finite codes. The proof is extended straightforward to rational codes.

\subsection{Bernouilli morphisms}

A Bernouilli morphism is a multiplicative morphism $\pi : A^* \rightarrow \mathbb{R}_+$ such that $\pi \, | \, A$ is a probability on $A$ such that $\pi(a) > 0$ for any $a$ in $A$.

\b

It is known that if $L$ is a language having the property that it does not intersect some ideal in $A^*$, then $\pi(L) = \sum_{w \in L} \, \pi(w) < \infty$.  This property is true if $L$ is rational code. See \cite{BP86} Prop. I.5.6 and Prop. I.5.12.

\b

From this, we deduce that $\pi(L) < \infty$ for each language $L = C_i, \, P_i, \, S, \, P, \, F$ considered in Section 4.2.

\subsection{Probabilistic free field}

We know that $\mathbb{Q}\langle A \rangle$ is embedded in the corresponding free field denoted $\mathcal{F}$.  Consider now the $\mathbb{Q}$-algebra $\mathbb{Q} \langle A \rangle / A - 1$, which is the quotient of $\mathbb{Q} \langle A \rangle$ by its two-sided ideal generated by $A - 1 = \sum_{a \in A} \, a - 1$. This $\mathbb{Q}$-algebra is a free associative algebra, since the relation $A = 1$ allows to eliminate one variable.  We denote it $\mathbb{Q} \langle A / (A - 1) \rangle$  Hence, there is a corresponding free field, denoted $\mathcal{P}$ and which we call the \textit{probabilistic free field}.

\begin{thm} \label{thm 2}
Let $\mu: A^* \rightarrow \mathbb{Q}^{n \times n}$ be the homomorphism corresponding to a complete and transitive unambiguous automaton.  Let $M = \sum_{a \in A} \, a \mu a$ be its matrix, $P$ the image of $M$ in the probabilistic free field $\mathcal{P}$, $C_i$ the code generating the fixpoints of vertex $i$, $P_i$ the sum of the labels of all paths starting at $i$ and not passing again through $i$.  Then $P_i, \, C_i$ and $P_i^{-1}$ can be evaluated in $\mathcal{P}$.  Moreover, the following equalities hold in $\mathcal{P}$:
\newcounter{marqueur5}
\begin{list}{(\roman{marqueur5})}{\usecounter{marqueur5}}
\item $C_i = 1$;
\item $(1 - t)(tP)^* \in \mathcal{P}(t)$ can be evaluated at $t = 1$ and its diagonal elements are $\lambda(C_i)^{-1}$, $i = 1, \, \hdots, \, n$.
\item $(P_1^{-1}, \, \hdots, \, P_n^{-1}) P = (P_1^{-1}, \, \hdots, \, P_n^{-1})$;
\item $\displaystyle{\sum_{i = 1}^n \, P_i^{-1} = 1}$;
\item for any maximal columns $\ell, \, \ell'$, $(P_1^{-1}, \, \hdots, \, P_n^{-1}) \, \ell = (P_1^{-1}, \, \hdots, \, P_n^{-1}) \, \ell'$.
\end{list}
\end{thm}

\n \textbf{Example 2.} (continued)

$C_1 = 1$ holds in $\mathcal{P}$, since one has, even in $\mathbb{Q} \langle A \rangle: C_1 - 1 = (1 + a) (a + b + c  -1) (1+b).$
Moreover, we have in $\mathcal{P}$
$$C_2 = (c + b^2 + cb)^*a + b(c + b^2 + cb)^*a = (1+ b)(c + b^2 + cb)^*a.$$
Now, in $\mathbb{Q} \langle\langle b, \, c\rangle\rangle$, one has $(1 + b)(c + b^2 + cb)^* = (b+c)^*$, since $\{c, \, b^2, \, cb\}$ is a complete suffix code with set of suffixes $\{1, \, b\}$ (see \cite{BP86}). Thus $C_2 = (b+c)^*a = 1$ since $a = 1 - b - c$.  Also,
\begin{eqnarray*}
C_3 & = & b(c + a)^*(b+ c+ a) \\
& = & b(c+a)^* = 1.
\end{eqnarray*}
In $\mathcal{P}$, we have $S = \begin{bmatrix}
c & a & b + c \\
1 & 0 & 1 \\
b & 0 & 0
\end{bmatrix}$.
We show that $P_2 = a^{-1}P_1$; indeed
\begin{eqnarray*}
P_2 & = & 1 + (c + b^2 + cb)^*(1 + b + c) + 1 + b(c+b^2 + cb)^*(b+c) + b(c + b^2 + cb)^* \\
& = & 2 + (1+b)(c + b^2 + cb)^*(1 + b + c) \\
& = & 2 + (b+c)^*(1 + b + c) \\
& = & 2 + (b+ c)^* + (b+c)^*(b+c) \\
& = & 1 + 2 (b + c)^* = 1 + 2 a^{-1} = a^{-1} P_1,
\end{eqnarray*}
since $P_1 = 2 + a$.  We deduce that $P_1^{-1} a = P_2^{-1}$.  Moreover
$$P_3 = 1 + b(c + a)^*(1+a) = 2 + a = P_1,$$
since $b(c+a)^* = 1$.  Thus
\begin{eqnarray*}
P_1^{-1}(b + c) + P_2^{-1} & = &  P_1^{-1}(a+ b+ c) = P_1^{-1} = P_3^{-1}, \\
P_1^{-1}c + P_2^{-1} + P_3^{-1}b & = & P_1^{-1}(c + a + b) = P_1^{-1}.
\end{eqnarray*}
This shows that $(P_1^{-1}, \, P_2^{-1}, \, P_3^{-1}) S = (P_1^{-1}, \, P_2^{-1}, \, P_3^{-1}) $.
Furthermore, $P_1^{-1} + P_2^{-1} + P_3^{-1} = P_1^{-1} (1 + a + 1) = 1$.  Now, the only two maximal columns are $\tgauche{(1, \, 1, \, 0)}$ and $\tgauche{(0, \, 1, \, 1)}$. We have
$(P_1^{-1}, \, P_2^{-1}, \, P_3^{-1}) \, \begin{bmatrix} 1 \\ 1 \\ 0 \end{bmatrix} = (P_1^{-1}, \, P_2^{-1}, \, P_3^{-1})   \, \begin{bmatrix} 0 \\ 1 \\ 1 \end{bmatrix} $ since $P_1 = P_3$.

\subsection{Proof of theorem}

We need the following lemma.

\begin{lem} \label{lemme: proprietes de P}
Let $P = \sum_{a \in A} \, a \mu a \in \mathcal{P}^{n \times n}$ be the image in $\mathcal{P}^{n \times n}$ of the matrix $M$ of some complete and transitive unambiguous automaton, with associated homomorphism $\mu: A^* \rightarrow \mathbb{Q}^{n \times n}$.  Then $P$ has the eigenvalue 1 with associated eigenspace of dimension 1.  Moreover, if $t$ is a central variable, then in $\mathcal{P}(t)$, $(1-t)(tP)^*$ can be evaluated for $t = 1$ and its rows span the eigenspace above.
\end{lem}

\n \textit{Proof}
Consider the (left) $\mathcal{P}$-subspace $E$ of $\mathcal{P}^{1 \times n}$ spanned by the maximal rows.  It has as subspace the subspace $E'$ spanned by the differences of such rows.  Let $C$ be the sum of the distincts columns of some element of the minimal ideal of $\mu A^*$. Then $rC = 1$ if $r$ is a maximal row (dual statement of Prop. 1).  Thus $E'$ is strictly included in $E$.

\b

By Section 4.1, for each maximal row $r$ and each $a \in A$, $r \mu a$ is a maximal row, denoted $r_a$.  Then
$$rP = \sum_{a \in A} \, r (\mu a) a = \sum_{a \in A} \, r_a a = r + \sum_{a \in A} \, (r_a - r) a,$$
since $\sum_{a \in A} a = 1$ in $\mathcal{P}$.  Thus $r$ is fixed by $P$ modulo the subspace $E'$.  Hence $P$ has 1 as eigenvalue.

\b

We show that its multiplicity is 1.  Indeed the multiplicity does not decrease under specialization.  For the latter, we take a positive Bernouilli morphism $\pi$; then $\pi(P)$ is an irreducible matrix because the automaton is transitive; it has nonnegative coefficients.  We claim that its eigenvalues are of module $\leq 1$.  Thus, we may apply the Perron-Frobenius theorem (\cite{LT} Section 15.3 Th.1) and, since 1 is an eigenvalue of $\pi(P)$ by the previous calculations, it is a root of multiplicity 1 of the characteristic polynomial.  But we know that 1 is an eigenvalue of $P$, hence it has multiplicity 1.  We conclude by using Lemma 1.

\b
 
It remains to prove the claim.  Since the automaton is unambiguous, the matrix $M^n = \l ( \sum_{a \in A} a \mu a \r )^n = \sum_{w \in A^n} \, w \mu w$ has the property that each entry is a subsum of $\sum_{w \in A^n} \, w$.  Hence each entry of $\pi(M^n) = \pi(P^n)$ is bounded by 1.  Hence each eigenvalue of $\pi(P)$ has module $\leq 1$.
\begin{flushright}
$\square$
\end{flushright}

\b

\n \textit{Proof of Theorem \ref{thm 2}}

$C_i$ is a rational maximal code.  So we may use the result at the end of Section 4.2: $A^* = S C_i^* P + F$, where $S, \, P, \, F$ are rational languages contained in the set of factors of $C_i$.  Then, by Section 4.3., $\pi(S), \, \pi(P)$ and $\pi(F)$ are $< \infty$ for any Bernouilli morphism.  This implies that $S, \, P, \, F$ can be evaluated in $\mathcal{P}$ (cf. the proof of Corollary 1).  The same holds for $C_i$ and $P_i$.  Now the equality in $\mathbb{Q}\langle\langle A \rangle\rangle$ above may be rewritten:
\begin{eqnarray*}
A^* - F & = & S C_i^* P \\
\R 1 - (1 - A) \, F & = & (1- A) \, S C_i^* P \\
\R (1 - (1 - A) \, F)^{-1} & = & P^{-1} (1 - C_i) S^{-1} A^* \\
\R 1 - C_i & = & P (1 - (1 - A) \, F)^{-1} (1 - A) S.
\end{eqnarray*}
This holds in $\mathbb{Q}\langle\langle A \rangle\rangle^{rat}$ and all these rational expressions can be evaluated in $\mathcal{P}$.  Thus, in $\mathcal{P}$, we obtain $1 - C_i = 0$.

\b

Let $Q_i$ be the set of paths from 1 to any vertex, that do not pass again through $i$.  Then we have, as in the proof of Th. \ref{thm: P_i sont definies dans un corps libre}, for any $i, \, j$,
$$(M^*)_{1i} \, P_i + Q_i = (M^*)_{1j} \, P_j + Q_j.$$
Arguing as in the latter proof, we find that,  denoting $\alpha_i$ the value of $(1-t)(tP)_{1i}^*$ at $t = 1$ (which exists by Lemma \ref{lemme: proprietes de P}), we obtain
$\alpha_i P_i = \alpha_j P_j$.
Note that $(tP)_{11}^* = C_1(t)^*$, where $C_1(t)^*$ denotes the canonical image of $\sum_{w \in C_1^*} \, t^{|w|} w \in \mathbb{Q}\langle A \rangle [[t]]^{rat}$ under the composition of homomorphisms
$$\mathbb{Q}\langle A \rangle [[t]]^{rat} \rightarrow \mathbb{Q}\langle A / (A - 1) \rangle  [[t]]^{rat} \rightarrow \mathcal{P} [[t]]^{rat} \rightarrow \mathcal{P}(t).$$
Thus $\alpha_1$ is the value at $t = 1$ of $(1 - t) C_1(t)^*$. Now taking the previous notations with $i = 1$, we have in $\mathcal{P}(t)$: $(tA)^* = S(t)C_1(t)^*P(t) + F(t)$.  Multiplying by $(1-t)$ and putting $t = 1$, we obtain, since $S, \, P, \, F$ can be evaluated in $\mathcal{P}$: $1 = S \alpha_1 P$.  Thus $\alpha_1 = S^{-1}P^{-1}$.

\b

Now, we have also $C_1 - 1 = P(1 - (1-A)F)^{-1} (A - 1)S$.  Thus, in $\mathcal{F}$, letting $P' = P(1 - (1-A)F)^{-1}$, $$ \lambda(C_1) = \lambda(P')(A - 1)S + P' \lambda(A) S + P'(A-1) \lambda(S).$$
We deduce that, in $\mathcal{P}$, $\lambda(C_1)  =P' \lambda(A) S = PS$.  This shows that $\alpha_1 = \lambda(C_1)^{-1}$.  This proves $(ii)$ and in particular, $\alpha_1 \neq 0$.  Thus, since $P_1 \neq 0$ in $\mathcal{P}$, all $\alpha_i$ and $P_i$ are $\neq 0$ in $\mathcal{P}$. Then $(iii)$ and $(iv)$ are proved as in the proof of Th. \ref{thm: P_i sont definies dans un corps libre}.

\b

In order to prove $(v)$, we observe that the elements of the minimal ideal $I$ of $\mu A^*$ are those of this monoid which have a minimal number of distinct nonnull rows (see \cite{BPR} Exercice VI.3.5 or \cite{BCKP} Proposition 1). This implies that if $r_1, \, \hdots, \, r_k$ are the distinct nonnull rows of some element $\mu w$ of $I$, then for any letter $a$, $r_1 \mu a, \, \hdots, \, r_k \mu a$ are the distinct nonnull rows of $\mu(wa)$.  We deduce that the span of the elements $r_1 + \hdots + r_k$ is invariant under the matrices $\mu a$.  Let $F$ denote this subspace, and $F'$ the subspace spanned by the difference of such elements.  By Prop.1, we have that $F'$ is strictly included in $F$.  Hence, there is a vector in $F$ fixed by each $\mu a$.  This implies that the eigenvector for eigenvalue 1 of the matrix $P$ is in $F$ and is therefore orthogonal to each difference of maximal columns of $\mu A^*$.  This proves $(v)$.
\begin{flushright}
$\square$
\end{flushright}

\section{Appendix 1: the commutative case}

The following result in an exercice on determinants.

\begin{lem} \label{lemme commutatif}
If the column eigenvector $\tgauche{(1, \, \hdots, \, 1)}$ is in the right kernel of a square matrix over a commutative ring, then the row vector $(m_1, \, \hdots, \, m_n)$ is in its left kernel, where $m_i$ is the $i$-th principal minor of the matrix.
\end{lem}

\b

From this, one may deduce the so-called \textit{Markov Chain tree theorem}, by using, as suggested in \cite{LW} page 4, the matrix-tree theorem, see e.g. \cite{S} Th. 5.6.8.

\b

The Markov chain tree theorem gives a formula, using spanning trees of the complete graph, for the stationary distribution of a finite Markov chain.  Equivalently, this formula gives a row vector fixed by a matrix fixing $\tgauche{(1, \, \hdots, \, 1)}$.  This theorem is attributed to Kirchoff by Persi Diaconis, who gives a probabilistic proof of it (see \cite{Bro89} p. 443 and 444).  See also \cite{AnTs89}, \cite{Ald90}.

\b

The Markov chain tree theorem is as follows: let $(a_{ij})$ be a stochastic matrix (that is, fixing $\tgauche{(1, \, \hdots, \, 1)}$).  Then the row vector $(b_1, \, \hdots, \, b_n)$ is fixed by this matrix, where $b_i$ is the sum of the weights of all spanning trees of the complete digraph on $\{1, \, \hdots, \, n\}$, rooted at $i$ (the edges of the tree all pointing toward $i$).  Here the weight of a subgraph is the product of the $a_{ij}$, for all edges $(i, \, j)$ in the subgraph.

\b

Using our Theorem \ref{thm: P_i sont definies dans un corps libre}, one easily deduces that $B = \sum_{i=1}^n \, b_i$, the sum of the weights of all rooted trees, is equal to the derivative of $\det(1 - (a_{ij}))$, with respect to the derivation fixing each $a_{ij}$.

\section{Appendix 2: Quasideterminants of stochastic matrices}

\subsection{Theory of quasideterminants}

The theory of quasideterminants was developed as a tool for linear algebra over noncommutative rings replacing the theory of determinants over commutative rings. Quasideterminants were introduced in \cite{Gel-Ret1} and developed in \cite{Gel-Ret2}, \cite{GR}, \cite{Gel-Ret4} and \cite{Gel-Ret5}. Let $R$ be an associative unital ring and let $A=(a_{ij})$, $i,j=1,2,\dots, n$ be a matrix over $R$. Denote by $A^{ij}$ the submatrix of $A$ obtained from from $A$ by deleting its $i$-th row and its $j$-th column. Set $r_i=(a_{i1},a_{i2},\dots , \hat a_{ij},\dots , a_{in})$ and $c_j=(a_{1j}, a_{2j},\dots , \hat a_{ij},\dots , a_{nj})$. Recall, that for any
matrix $C$ we denote by ${}^tC$ the transposed matrix.

\b

\n \textbf{Definition}. Suppose that the matrix $A^{ij}$ is invertible. Then the quasideterminant $|A|_{ij}$ is defined as
$$|A|_{ij}=a_{ij} - r_i(A^{ij})^{-1}c_j.$$

\n \textbf{Example}. If $n=2$, then $|A|_{12}=a_{12}-a_{11}a_{21}^{-1}a_{22}.$

\b

Let matrix $A$ be invertible and $A^{-1}=(b_{pq})$. If the quasideterminant $|A|_{ij}$ is invertible then $b_{ji}=|A|_{ij}^{-1}.$ In the commutative case, $|A|_{ij}=(-1)^{i+j}\det A/\det A^{ij}$.

\b

It is sometimes convenient to use another notation for quasideterminants $|A|_{ij}$ by boxing the leading entry, i.e.
$$ |A|_{ij}=
\begin{vmatrix}
\dots &\dots &\dots \\
\dots &\boxed{a_{ij}}&\dots \\
\dots &\dots &\dots \end{vmatrix}.$$

We remind now the basic properties of quasideterminants (see \cite{Gel-Ret1}). An equality $|A|_{pq}=|B|_{rs}$ means that the first quasideterminant is defined if and only if the second quasideterminant is defined and that both quasideterminants
are equal. The properties are:

\b

\newcounter{marqueur15}
\begin{list}{(\roman{marqueur15})}{\usecounter{marqueur15}}
\item \textit{Permutations of rows and columns}:
Let $\sigma , \tau$ be permutations of $\{1,2,\dots ,n\}$. Set $B=(a_{\sigma (i), \tau (j)})$. Then  $|A|_{pq}=B_{\sigma (p),\tau (q)}$. \\
\item \textit{Multiplication of row and columns}:

\b

Let the matrix $B=(b_{ij})$ be obtained from matrix $A$ by
multiplying the $i$-th row by $\lambda \in R$ from the left, i.e. $b_{ij}=\lambda a_{ij}$ and $b_{kj}=a_{kj}$ for all $j$ and $k\neq i$. Then $|B|_{kj}=\lambda |A|_{kj}$ if $k=i$, and $|B|_{kj}=|A|_{kj}$ if $k\neq i$ and $\lambda $ is invertible.

\b

Let the matrix $C=(c_{ij})$ be obtained from matrix $A$ by multiplying the $j$-th column by $\mu \in R$ from the right, i.e. $c_{ij}=a_{ij}\mu $ and $c_{il}=a_{il}$ for all $i$ and $l\neq j$. Then $|C|_{il}=|A|_{il}\mu$ if $l=j$, and
$|C|_{il}=A|_{il}$ if $l\neq j$ and $\mu $ is invertible. \\

\item \textit{Addition of rows and columns}:

Let the matrix $B$ be obtained from $A$ by replacing the $k$-th row of $A$ with the sum of $k$-th and $l$-th row, i.e. $b_{kj}=a_{kj}+a_{lj}$, $b_{ij}=a_{ij}$ for $i\neq k$. Then
$|A|_{ij}=|B|_{ij}$, $i=1,2,\dots , \hat k,\dots , n$, $j=1,2,\dots ,n$.

\b

Let the matrix $C$ be obtained from $A$ by replacing the $k$-th column of $A$ with the sum of $k$-th and $l$-th column, i.e. $c_{ik}=a_{ik}+a_{il}$, $b_{ij}=a_{ij}$, $c_{ij}=a_{ij}$ for $j\neq k$. Then $|A|_{ij}=|C|_{ij}$, $i=1,2,\dots , n$, $j=1,2,\dots ,\hat l, \dots n$.
\end{list}

\b

In \cite{Gel-Ret1} a noncommutative analogue of the Cramer's rule for systems of left linear equations, i.e. the systems when coefficients are at the left of the unknowns, was formulated. The analogue for systems of right linear equations can be formulated as follows.

\b

Let $B=(b_{ij})$ an $n\times n$-matrix over $R$, $\xi=(\xi_i)$ be a row-matrix over $R$ and $x=(x_i)$ be a row-matrix of unknowns. Here $i,j=1,2,\dots , n$. For $1\leq k\leq n$ denote by $B(\xi, k)$ the matrix obtained from $B$ by replacing the $k$-th row of $B$ by $\xi $.

\begin{prop}
If $xB=\xi $ then $$x_k|B|_{kq}=|B(\xi,k)|_{kq}$$
for any $k$ provided that the both quasideterminants are defined.
\end{prop}

\b

\n \textbf{Example}. For $n=2$ one has
$$x_1(b_{12}-b_{11}b_{21}^{-1}b_{22})=\xi_2 -\xi_1b_{21}^{-1}b_{22}$$ and
also
$$x_1(b_{11}-b_{12}b_{22}^{-1}b_{21})=\xi_1 -\xi_1b_{22}^{-1}b_{21}.$$

\subsection{Results}
\begin{lem} \label{lemme 2 de Retakh}
Let $A=(a_{ij})$, $i,j=1,2,\dots ,n$ be a stochastic matrix over $R$.
Consider the system of $n+1$ equations
\begin{equation}
\sum _{i=1}^n \, x_i \, a_{ij}=x_j,
\label{sum i=1..n x_i a_(ij) = x_j}
j=1,2,\dots , n
\end{equation}
together the equation
$$\sum _{i=1}^nx_i=1.$$
Then any of the n first equations of the system is a corollary of the other $n$ equations.
\end{lem}

\b

\n \textit{Proof} Take any $1\leq k\leq n$ and add all equations (\ref{sum i=1..n x_i a_(ij) = x_j}) for $j\neq k$. The right hand side of the sum can be written as $1-x_k$ and the left hand side can be written as $\sum
_{i=1}^nx_i(1-a_{ik})$. As a result we have
$$\sum _{i=1}^nx_i(1-a_{ik})=1-x_k$$
which implies
$$\sum _{i=1}^nx_ia_{ik}=x_k.$$
The lemma is proved.
\begin{flushright}
$\square$
\end{flushright}

\begin{thm}\label{thm 3 de Retakh}
Let $A=(a_{ij})$, $i,j=1,2,\dots , n$  be a stochastic matrix. The system
$$\sum _{i=1}^nx_ia_{ij}=x_j, \ \ \ j=1,2,\dots , n,$$
$$\sum _{i=1}^nx_i=1$$
has a unique solution over the algebra of series in variables $a_{ij}$ satisfying the relations $\sum _{j=1}^na_{ij}=1$.
The solutions are given by the formula
$$x_k^{-1}=1+\sum a_{ki_1}a_{i_1i_2}a_{i_2i_3}\dots a_{i_{s-1}i_s}$$
where the sum is taken over all sets of naturals $i_1,i_2,\dots , i_s$ where $s\geq 1$ and $i_p\neq k$, $p=1,2,\dots , s$.
\end{thm}

\n \textit{Proof} Lemma \ref{lemme 2 de Retakh} implies that $x_1,\dots , x_n$ are solutions of the system
$$\sum_{i=1}^n x_i = 1, $$
$$\sum_{i=1}^n x_i (a_{ij}-\delta _{ij}) = 0, \q j \neq k.$$

Write the system in the form $xB=\xi$, where $x=(x_1,\dots , x_n)$ and $\xi=(1,0,\dots ,0)$, and apply Proposition 11.  Note that $|B(\xi , k)|_{k1}=1$: indeed, the $k$-th row of $B(\xi , k)$ is $\xi$ and therefore, by the definition of quasideterminants, $r_k=(0,\dots,0)$ and $|B(\xi , k)|_{k1}=B(\xi , k)_{k1}=1$. Therefore,
$$x_k|B|_{k1}=1.$$

\b

Recall that $A^{kk}$ is the submatrix of $A$ obtained from $A$ by omitting its $k$-th row and $k$-th column. Set $C=A^{kk}$. Let $I$ be the unit matrix of order $n-1$ and $a(k)=(a_{k1}, a_{k2},\dots , \hat a_{kk},\dots , a_{kn})$. Note that the first column of matrix $B$ is $\gamma={}^t(1,1,\dots , 1)$. By the definition of quasideterminants
$$|B|_{k1}=1 - a(k)(C-I)^{-1}\gamma=1 + \sum _{p=0}^{\infty}C^p=1+\sum a_{ki_1}a_{i_1i_2}a_{i_2i_3}\dots a_{i_{s-1}i_s}.$$
This proves the theorem.
\begin{flushright}
$\square$
\end{flushright}

\b
\begin{rem}
Note that monomials $a_{ki_1}a_{i_1i_2}\dots a_{i_{s-1}i_s}$ can be interpreted as paths in the complete graph with vertices
$1,2,\dots ,n$.
\end{rem}
\bigskip

\subsection{Stochastic matrices and main quasiminors}

\b

Observe that a matrix $M$ is stochastic if and only if $\tgauche{(1, \, 1, \, \hdots, \,1)}$ is in the kernel of $M-I$. This justifies the next results.

\b

\begin{lem}\label{lemme:matrice stochastique sur un anneau unitaire}
Let $A=(a_{ij})$, $i,j=1,2,\dots , n$ over an associative unital ring annihilate the column vector $\tgauche{(1, \, 1, \, \hdots, \,1)}$. For $p\neq q$ one has
$$|A^{pq}|_{qp}=-|A^{pp}|_{qq}$$
if the right hand side is defined.
\end{lem}

\b

\n \textit{Proof} Without loss of generality one can assume that $p=1$ and $q=n$. Then
$$ |A^{1n}|_{n1}=\begin{vmatrix}
a_{21} &a_{22}&\dots &a_{2n-1} \\
a_{31} &a_{32} &\dots &a_{3n-1}\\
  & &\dots & & \\
\boxed{a_{n1}} &a_{n2}&\dots &a_{nn-1}\end{vmatrix}.$$

Since $A$ is stochastic by hypothesis, we can rewrite the elements $a_{k1}$, $k=2,3,\dots , n$ as $a_{k1}=-a_{k2}-a_{k3}-\dots -a_{kn}$. By adding the columns in the last quasideterminant to the first one and using property (iii), property (ii) for the first column and $\mu =-1$ and property (i) we  get the expression
$$
-\begin{vmatrix}
a_{2n} &a_{22}&\dots &a_{2n-1} \\
a_{3n} &a_{32} &\dots &a_{3n-1}\\
  & &\dots & & \\
\boxed{a_{nn}} &a_{n2}&\dots &a_{nn-1}\end{vmatrix}=
-\begin{vmatrix}
a_{22} &a_{23}&\dots &a_{2n} \\
a_{32} &a_{33} &\dots &a_{3n}\\
  & &\dots & & \\
a_{n2} &a_{n3}&\dots &\boxed{a_{nn}}\end{vmatrix}$$
which is $-|A^{11}|_{nn}.$ Our computations also show the existence of $|A^{1n}|_{n1}$. The lemma is proved.
\begin{flushright}
$\square$
\end{flushright}

\b

\begin{thm} \label{thm 5 de Retakh}
Let $A=(a_{ij})$, $i,j=1,2,\dots , n$ annihilate the column vector $\tgauche{(1, \, 1, \, \hdots, \,1)}$. Assume that all quasideterminants $|A^{ii}|_{jj}=$ are defined for $i\neq j$. Then $xA=0$ where $x=(x_1,x_2,\dots , x_n)$ if and only if
$$x_i|A^{jj}|_{ii}=x_j|A^{ii}|_{jj}, \ \ i\neq j.$$
\end{thm}

\b

\begin{rem}
In the commutative case $|A_{ii}|_{jj}=\det A^{ii}/\det A^{ij,ij}$ provided that the denominator is defined.
Here $A^{ij,ij}$ is the submatrix of $A$ obtained from $A$ by removing its rows and columns with the indices $i$ and $j$. Thus the theorem implies
$$x_im_j=x_jm_i$$
where $m_i$'s are the main minors of the matrix and we may choose $x_i=m_i$, $i=1,2,\dots ,n$ as a solution of the equation $xA=0$ obtaining Lemma \ref{lemme commutatif}.
\end{rem}

\b

\n \textit{Proof of the theorem}. We will prove the "if" part. The "only if" part can be proved by reversing the arguments. Without loss of generality, we assume that $i=1$ and $j=n$. Note that $x_i$'s satisfy the system of linear equations
$$ \sum _{p=1}^{n-1}x_pa_{pq}=-x_na_{nq},\ \ \ q=1,2,\dots , n-1. $$

The Cramer's rules give us the equality
$$ x_1|A^{nn}|_{11}=\begin{vmatrix}
\boxed{-x_na_{n1}} & -x_na_{n2} &\dots & -x_na_{nn-1}\\
a_{21} & a_{22}&\dots & a_{2n-1} \\
  & &\dots & &  \\
a_{n-11} &a_{n-12}&\dots &a_{n-1n-1}  \end{vmatrix}.$$

By properties (i) and (ii), the right hand side equals to $-x_n|A^{1n}|_{n1}$.

\b

The theorem now follows from Lemma \ref{lemme:matrice stochastique sur un anneau unitaire}.
\begin{flushright}
$\square$
\end{flushright}

By using the results from \cite{Gel-Ret2} and \cite{GR} we can show that Theorem \ref{thm 5 de Retakh} implies Theorem \ref{thm 3 de Retakh} provided the corresponding quasideterminants are invertible.

\end{document}